\documentclass[11pt]{article}

\usepackage{amsfonts,epsfig,fullpage}
\usepackage{caption}
\usepackage{amsmath}
\usepackage[table]{xcolor}
\usepackage{multirow}
\usepackage{graphics}

\newtheorem{theo}{Theorem}

\newtheorem{lem}{Lemma}
\newtheorem{claim}[lem]{Claim}
\newtheorem{coro}[lem]{Corollary}

\newtheorem{define}{Definition}

\newtheorem{alg}{Algorithm}
\newtheorem{proc}{Procedure}

\newcommand{\BE}{\begin{enumerate}} \newcommand{\EE}{\end{enumerate}}
\newcommand{\BI}{\begin{itemize}} \newcommand{\EI}{\end{itemize}}
\newcommand{\BDes}{\begin{description}}\newcommand{\EDes}{\end{description}
}
\newcommand{\BT}{\begin{theo}} \newcommand{\ET}{\end{theo}}
\newcommand{\BL}{\begin{lem}} \newcommand{\EL}{\end{lem}}
\newcommand{\BD}{\begin{define}} \newcommand{\ED}{\end{define}}
\newcommand{\BCM}{\begin{claim}} \newcommand{\ECM}{\end{claim}}
\newcommand{\BC}{\begin{coro}} \newcommand{\EC}{\end{coro}}
\newcommand{\BA}{\begin{alg}} \newcommand{\EA}{\end{alg}}
\newcommand{\BP}{\begin{proc}} \newcommand{\EP}{\end{proc}}

\def\FullBox{\hbox{\vrule width 8pt height 8pt depth 0pt}}
\newcommand{\qed}{\;\;\;\FullBox}
\newenvironment{proof}{\noindent{\bf Proof:~~}}{\(\qed\)}
\newcommand{\BPF}{\begin{proof}} \newcommand {\EPF}{\end{proof}}
\newenvironment{proofof}[1]{\noindent{\bf Proof of {#1}:~~}}{\(\qed\)}
\newcommand{\BPFOF}{\begin{proofof}} \newcommand {\EPFOF}{\end{proofof}}

\newcommand{\BEQ}{\begin{equation}} \newcommand{\EEQ}{\end{equation}}
\newcommand{\BEQN}{\begin{eqnarray}}\newcommand{\EEQN}{\end{eqnarray}}

\newlength{\saveparindent}
\newlength{\saveparskip}
\newcommand{\eqdef}{\stackrel{\rm def}{=}}

\providecommand{\keywords}[1]{\textbf{\textit{Key Words---}} #1}
\providecommand{\subjclass}[1]{\textbf{\textit{Classification Codes---}} #1}

\begin{document}

\title{On maximal isolation sets in the uniform intersection matrix}

\author{
    Michal Parnas \\
    The Academic College\\
    of Tel-Aviv-Yaffo \\
     {\tt michalp@mta.ac.il}
\and
    Adi Shraibman \\
   The Academic College\\
    of Tel-Aviv-Yaffo \\
     {\tt adish@mta.ac.il}
}

\maketitle

\abstract{
Let $A_{k,t}$ be the matrix that represents the adjacency matrix of the intersection bipartite graph of
all subsets of size $t$ of $\{1,2,...,k\}$. We give constructions of large isolation sets in $A_{k,t}$,
where, for a large enough $k$, our constructions are the best possible.

We first prove that the largest identity submatrix in $A_{k,t}$ is of size $k-2t+2$.
Then we provide constructions of isolations sets in $A_{k,t}$ for any $t\geq 2$, as follows:
\begin{itemize}
\item
If $k = 2t+r$ and  $0 \leq r \leq 2t-3$, there exists an isolation set of size $2r+3 = 2k-4t+3$.
\item
If $k \geq 4t-3$, there exists an isolation set  of size $k$.
\end{itemize}
The construction is maximal  for $k\geq 4t-3$, since the Boolean rank of $A_{k,t}$ is $k$ in this case.
As we prove, the construction is maximal also for  $k = 2t, 2t+1$.

Finally, we consider the problem of the maximal triangular isolation submatrix of $A_{k,t}$
that has ones in every entry on the main diagonal and below it, and zeros elsewhere.
We give an optimal construction of such a submatrix of size $({2t \choose t}-1) \times ({2t \choose t}-1)$,
for any $t \geq 1$ and a large enough $k$.
This construction is tight, as there is a matching upper bound, which can be derived from a theorem of Frankl about skew matrices.
}\\

\keywords{Boolean rank; isolation set; intersection matrix.}

\subjclass{15A23, 15B34, 15A03}

\section{Introduction}
Let $A_{k,t}$ be the matrix that represents the adjacency matrix of the intersection bipartite graph of
all subsets of size $t$ of $[k] \eqdef \{1,2,...,k\}$. Thus, each row and column of $A_{k,t}$ is
indexed by a subset in ${[k] \choose t}$.
The size of $A_{k,t}$ is ${k \choose t}\times {k \choose t}$,
and $A[x][y] = 1$ if and only if the two subsets $x,y$ intersect.

Intersecting families of subsets have been studied extensively over the years,
and some of the results achieved can be inferred as results about families of submatrices of the matrix $A_{k,t}$.
For example, Pyber~\cite{pyber1986new} proved that the maximal cross-intersecting family of subsets of ${[k] \choose t}$
is of size ${k-1 \choose t-1}^2$, and thus, the largest all-ones submatrix of $A_{k,t}$ is of size ${k-1 \choose t-1}^2$.
Another example is a theorem of Bollob{\'a}s~\cite{bollobas1965generalized} about cross-intersecting sets,
that allows to show that the largest submatrix representing a crown graph in $A_{k,t}$ is of size ${2t \choose t}\times {2t \choose t}$.

Here we suggest to continue and explore various families of maximal submatrices of $A_{k,t}$.
In particular, we would like to find small submatrices of $A_{k,t}$ whose Boolean rank is large.
The {\em Boolean rank} of a matrix $B$ of size $n \times m$ is equal to the smallest integer $r$,
such that $B$ can be factorized as a product of two Boolean matrices, $XY=B$,   where
$X$ is a matrix of size $n \times r$  and $Y$ is a matrix of size $r \times m$, and all additions and multiplications are Boolean.
The Boolean rank is also equal to the  minimal number of monochromatic combinatorial rectangles
required to cover all of the ones of $B$, and it is equal to the minimal number of complete bi-cliques
needed to cover the edges of the bipartite graph whose adjacency matrix is $B$ (see~\cite{Gregory}).
Lastly, the Boolean rank is also tightly related to the notion of nondeterministic communication complexity \cite{KN97}.

The Boolean rank of $A_{k,t}$ was shown in~\cite{parnas2019boolean} to be $k$ for any $1 \leq t \leq k/2$.
Furthermore, it was proved in~\cite{parnas2019boolean} that there exists a family of submatrices of $A_{k,t}$,
each of size $(m \cdot s)\times (m \cdot s)$, where $m = {2t-2 \choose t-1}$ and $s= k-2t+2$,
whose Boolean rank is also $k$, for a large range of values of $k,t$.
These submatrices are rather large, and a question that arises is
if there are smaller submatrices of $A_{k,t}$ whose Boolean rank is $k$, or as close as possible to $k$.
We answer this question and prove that for a large enough $k$, there are, in fact, submatrices
of size $k\times k$ of $A_{k,t}$, whose Boolean rank is $k$.

Natural candidates for small matrices with a large Boolean rank are {\em isolation sets}
(or {\em fooling sets} as they are called in communication complexity).
An isolation set for a Boolean matrix $B$ is a subset of entries $F$ in $B$ that are all ones of $B$,
such that no two ones in $F$ are in the same row or column of $B$, and no two ones in $F$ are contained
in an all-one submatrix of size $2\times 2$ of $B$.
Throughout the paper we will represent an isolation set of a given matrix $B$ as a submatrix $F$ of $B$,
where the ones of the isolation set are on the main diagonal of $F$, and $F$ is called an {\em isolation matrix}.
The Boolean rank of an isolation matrix of size $f\times f$ is equal to $f$, and therefore, the size of the maximal isolation
set in a given matrix, bounds below the Boolean rank of that matrix (see for example~\cite{KN97, beasley2012isolation}).
Hence, finding large isolation sets in $A_{k,t}$ answers partially the question of finding small submatrices of $A_{k,t}$ with a large Boolean rank.

If $k < 2t$, then $A_{k,t}$ is just the all-ones matrix, since every two subsets of size $t$ intersect,
and thus, the largest isolation set is of size $1$. Therefore, the question of
finding large isolation sets in $A_{k,t}$ is interesting only for $k \geq 2t$.
The simplest form of an isolation matrix is the identity matrix, and thus, we first consider
the problem of determining the size of the largest identity submatrix in $A_{k,t}$, and prove the following:

\BT
\label{theorem1}
The largest identity submatrix in $A_{k,t}$ is of size $s\times s$, where $s = k-2t+2$.
\ET

Recall that the complement of $A_{k,t}$ is the adjacency matrix of the Kneser graph $KG_{k,t}$, in which the vertices are all
subsets of size $t$ of $[k]$, and there is an edge between two subsets $x,y$ if and only if $x \cap y = \emptyset$.
Furthermore, the complement of the identity matrix is the adjacency matrix of the crown graph of the same size.
Thus, from Theorem~\ref{theorem1},
we immediately get that the largest submatrix representing a crown graph in $KG_{k,t}$ is of size $s = k-2t+2$.
In particular, this is also the maximal size of a clique in $KG_{k,t}$,
which corresponds to the fact that the chromatic number of $KG_{k,t}$ is $k-2t+2$ \cite{lovasz1978kneser}.

Another simple isolation matrix is the {\em triangular matrix} with ones in every entry on the main diagonal and below it,
and zeros elsewhere.
We give an optimal construction of such a triangular matrix in $A_{k,t}$, where our construction  uses similar ideas to those used by
Tuza~\cite{tuza1987inequalities}, and the upper bound follows easily from a result of Frankl~\cite{frankl1982extremal}
that proved a skew version of a theorem of Bollob{\'a}s~\cite{bollobas1965generalized}.

\BT
\label{theorem2:triangular}
For any $t \geq 1$ and a large enough $k$, the maximal triangular submatrix of $A_{k,t}$ is of size
$d \times d$, where $d = {2t \choose t} -1$.
\ET

As can be seen, the size of the maximal triangular submatrix of $A_{k,t}$ does not depend on $k$ (as long as $k$ is large enough).
Thus, for a large enough $k$, the maximal identity submatrix $I_s$ promised by Theorem~\ref{theorem1},
is a larger isolation submatrix in $A_{k,t}$. But is $I_s$ the largest isolation matrix in $A_{k,t}$?
If $t=1$ then $A_{k,t} = I_k = I_s$, and in this case, this is, of course, the maximal isolation set.
It is also not hard to verify that if $k = 2t$, there exists an isolation set of size $k-2t+3 = 3$ in $A_{k,t}$,
and this is the maximal isolation set in this case (see for example~\cite{beasley2012isolation}).
As we prove, for $2 \leq t < k/2$, there are larger isolation sets,
and the submatrix $I_s$ is not the largest isolation matrix for these values of $t$ and $k$.
In fact, when $k$ is large enough, there exists in $A_{k,t}$ an isolation set of size $k$.

\BT
\label{theorem:constructions}
For any $t\geq 2$, the matrix $A_{k,t}$  has an isolation set of the following size:
\begin{itemize}
\item
If $k = 2t+r$ and  $0 \leq r \leq 2t-3$, there exists an isolation set of size $2r+3 = 2k-4t+3$.
\item
If $k \geq 4t-3$  there exists an isolation set of size $k$.
\end{itemize}
\ET

Notice that for any fixed given $t$, the size of the isolation set starts at $3$ when $k = 2t$,
and then grows by an additive term of two when $k$ is increased by one, until the point that $k = 4t-3$.
Then, we get an isolation set of maximal size $k$. Our construction is also maximal for $k = 2t$, and as we prove it
is also maximal for $k = 2t+1$. It is an open question if the construction is  maximal for $2t+2 \leq k \leq 4t-4$.

\BT
\label{theorem:upperbounds}
If $k = 2t+1$ and $t\geq 2$, then the size of any isolation set in  $A_{k,t}$ is at most $5$.
\ET

\section{The maximal identity submatrix in $A_{k,t}$}

In all that follows we denote the identity matrix of size $n\times n$ by $I_n$,
and refer to the subsets representing a row or column of $A_{k,t}$ as row or column indices.
Therefore, each row or column index is a subset of ${[k] \choose t}$.
We now prove Theorem~\ref{theorem1},
and show that the maximal identity submatrix  of $A_{k,t}$ is of size $s \times s$, where $s = k-2t+2$.

First notice that there exists such a large identity submatrix in $A_{k,t}$.
Just take $s$  row indices of the form $\{1,2,...,t-1\}\cup \{i\}$
and column indices of the form $\{t,t+1,...,2t-2\}\cup \{i\}$, for $i = 2t-1,2t,...,k$.
This defines an identity submatrix of $A_{k,t}$ of size $s\times s$.

We next show that this is the largest identity submatrix possible in $A_{k,t}$.
Clearly this is true for a submatrix on the main diagonal of $A_{k,t}$.
Assume, by contradiction, that there exists an identity submatrix $I_{s+1}$  on the main diagonal of $A_{k,t}$,
and let $x_1,...,x_{s+1}$ be the row and column indices of $I_{s+1}$,
where we have that $x_i \cap x_j = \emptyset$ if and only if $i \neq j$.
But then we get an independent set of size $s+1$ in $A_{k,t}$ that includes $x_1,...,x_{s+1}$.
Thus, the complement of $A_{k,t}$, that is, the Kneser graph $KG_{k,t}$, has a clique of size $s+1$.
This is in contradiction to the fact that the chromatic number of $KG_{k,t}$ is $s$ (see~\cite{lovasz1978kneser}).
In general though, the identity submatrix does not have to be on the main diagonal of $A_{k,t}$,
and thus, a different proof is needed.

We first need the following claim proved in~\cite{prs} that characterizes the decompositions of the identity matrix.
For completeness we include its proof.

\BCM[\cite{prs}]
\label{clm:identity}
Let $XY = I_n$ be a Boolean decomposition of the $n\times n$ identity matrix $I_n$, where
 $X$ is an $n \times r$ Boolean matrix and $Y$ is an $r \times n$ Boolean matrix. Denote by
$x_1, \ldots, x_r$ the columns of $X$ and by $y_1,\ldots,y_r$ the rows of $Y$.
Then:
\begin{enumerate}
\item
For each $i \in [r]$, either $x_i = y_i = e_j$ for some $j \in [n]$,
where $e_j$ denotes the $j^{\rm th}$ standard basis vector,
or $x_i$ is the all-zeros vector, or $y_i$ is the all-zeros vector.
\item
Furthermore, for each $j \in [n]$, there exists some $i \in [r]$ such that $x_i = y_i = e_j$.
\end{enumerate}
\ECM
\BPF
If we write the decomposition $XY = I_n$ with outer products, then $I_n = \sum_{i=1}^r x_i \otimes y_i$,
where $x \otimes y$ denotes the outer product of a column vector $x$ and a row vector $y$, i.e. it
is a matrix of size  $n\times n$.

Assume first that there exists an index $\ell \in [r]$ for which Item $1$ of the claim does not hold.
But then the matrix $x_{\ell} \otimes y_{\ell}$ contains a one that is not on the main diagonal of the matrix,
and since the addition is the Boolean addition, the sum $\sum_{i=1}^r x_i \otimes y_i \neq I_n$. Thus, Item 1 always holds
for any decomposition $XY$ of $I_n$.

Now assume that there exists some $j \in [n]$, such that there is no $i \in [r]$ for which $x_i = y_i = e_j$.
But then the $j^{\rm th}$ entry on the main diagonal of $\sum_{i=1}^r x_i \otimes y_i$ will be a zero.
\EPF\\

\BCM
\label{claim:numofones}
Let $XY = I_n$ be a decomposition of the $n \times n$ identity matrix $I_n$, where $X$ is an $n \times r$ Boolean matrix
and $Y$ is an $r \times n$ Boolean matrix. Then the total number of $1$'s in both $X$ and $Y$ is at most $2n + (r-n)n$.
\ECM
\BPF
By Claim~\ref{clm:identity}, for each $j \in [n]$, there exists some $i \in [r]$ such that $x_i = y_i = e_j$.
Assume, without loss of generality, that $x_i = y_i = e_i$ for $i = 1,...,n$.
Then the maximal number of $1$'s in any decomposition of $I_n$ occurs when
for all the remaining indices, $n < i \leq r$, it holds that
one of $x_i$ or $y_i$ is the all-zero vector and the other is the all-one vector.
Therefore, the number of ones in both $X$ and $Y$ is at most $2n + (r-n)n$.
\EPF\\

\BL
\label{claim:identity_size}
The largest identity submatrix of $A_{k,t}$ is $I_s$, where $s= k-2t+2$.
\EL
\BPF
Let $I_{\ell}$ be any identity submatrix of $A_{k,t}$.
Consider now the decomposition $XY = A_{k,t}$  of $A_{k,t}$, where $X = {[k] \choose t}$ and $Y = {[k] \choose t}$,
and let $X' \subseteq X, Y' \subseteq Y$, such that $X'Y' = I_{\ell}$.

Notice that $X'$ is an $\ell \times k$ matrix and $Y'$ is an $k \times \ell$ matrix, and  the total number
of $1$'s in both $X'$ and $Y'$ is exactly $2\ell t$.
But, by Claim~\ref{claim:numofones}, the total number of $1$'s in both $X'$ and $Y'$ is at most $2\ell+(k-\ell)\ell$. Thus,
$ 2\ell t \le 2\ell+(k-\ell)\ell$, and therefore $\ell \le k-2t+2$ as claimed.
\EPF\\

The following bound on the largest crown graph that is a submatrix of $KG_{k,t}$ is an immediate consequence
of Lemma~\ref{claim:identity_size}.
\begin{coro}
The largest matrix representing a crown graph that is a submatrix of the Kneser matrix $KG_{k,t}$,
 is of size $s\times s$, where $s = k-2t+2$.
\end{coro}

\section{Maximal triangular matrices in $A_{k,t}$}

As stated in the introduction, the following theorem of  Bollob{\'a}s~\cite{bollobas1965generalized}, allows to show that the largest
submatrix representing a crown graph in $A_{k,t}$, is of size ${2t \choose t}\times {2t \choose t}$, and this result is tight.
That is, there exists a simple construction of such  a large submatrix in $A_{k,t}$.

\BT[\cite{bollobas1965generalized}]
\label{bollobas}
Let $(A_i,B_i)$ be pairs of sets, such that $|A_i| = a, |B_i| =b$ for $1 \leq i \leq m$,
and  $A_i \cap B_j = \emptyset$ if and only if $i=j$. Then
$$m  \leq {a+b \choose a}.$$
\ET

This theorem has several generalizations, among them is a result of
Frankl~\cite{frankl1982extremal} that considered the skew version of the problem,
and showed that the same bound holds even under the following relaxed assumptions:
Let $(A_i,B_i)$ be pairs of sets, such that $|A_i| = a, |B_i| =b$ for $1 \leq i \leq m$,
$(A_i,B_i) = \emptyset$ for every $1 \leq i \leq m$, and $A_i \cap B_j \neq \emptyset$ if  $i>j$.
Then $m  \leq {a+b \choose a}$.
Note that for this formulation of the problem, all entries below the main diagonal
are ones, but above the main diagonal there can be either zeros or ones.

Here we consider the following special case:
What is the maximal number $m$ of pairs of subsets $(A_i,B_i)$,
such that $|A_i| = |B_i| = t$  for every $1 \leq i \leq m$,
and $A_i \cap B_j \neq \emptyset$ if and only if $i \geq j$.
Such a set of pairs of subsets defines a {\em triangular} submatrix of $A_{k,t}$ of size $m \times m$, for some large enough $k$,
with ones on the main diagonal and below it, and zeros elsewhere.
Denote such a matrix by $D_{m}$, and notice that $D_{m}$ is an isolation matrix.

Using the result of Frankl~\cite{frankl1982extremal} stated above,
it can be shown that the size of any triangular submatrix of $A_{k,t}$ is bounded above by $d \times d$, where $d ={2t \choose t}-1$.
To verify this, simply add to any maximal triangular submatrix an  additional first row and  last column that are all zero
(for a large enough $k$, it is always possible to define one more row index and column index that do not intersect with any of the given row and
column indices of the submatrix).
Thus, we get a matrix in which the main diagonal is all-zero, and below the main diagonal all elements are one.
By the result of Frankl, the size of such a matrix is at most  ${2t \choose t}\times {2t \choose t}$.
Hence, the size of any maximal triangular matrix is bounded above by  $({2t \choose t}-1)\times ({2t \choose t}-1)$.

We now proceed to prove Theorem~\ref{theorem2:triangular},
and show a construction of a triangular submatrix of $A_{k,t}$ that matches the above upper  bound.
The construction we describe is recursive, using an idea similar to what was done
by Tuza~\cite{tuza1987inequalities}.

Let $f(a,b)$ be the maximal $m$, such that $D_m$ is a submatrix of $A_{k,t}$ for a large enough $k$, having
row indices that are subsets of size $a$ and column indices that are subsets of size $b$.
We want to find $f(t,t)$, for $t\geq 1$.
We first give the stopping conditions for the recursion for $f(a,b)$.

\BCM
\label{f_1}
For any $a\geq 1$  it holds that $f(a,1) = f(1,a) = a$.
\ECM
\BPF
To see that $f(a,1) \geq a$, take as row indices the subsets $\{1, a+1,\ldots,2a-1\},\{1,2,a+1,\ldots,2a-2\},...\{1,2,...,a\}$,
and as column indices the subsets $\{1\},\{2\},...,\{a\}$.

For the lower bound, assume by contradiction that $f(a,1) = a+1$, and let the column indices be $\{1\},\{2\},...,\{a+1\}$.
Since the last row of the matrix is all-ones, then the index of the last row intersects with
all column indices. Thus, it contains the subset $\{1,2,..., a+1\}$,
in contradiction to the fact that the size of the row indices is $a$.
Hence, $f(a,1)=a$.

Similar arguments hold for $f(1,a)$, while exchanging the row and column indices.
\EPF\\

Now we can prove the general recursive formula for $f(a,b)$.
\BL
\label{const:triangular}
For any $a,b \geq 1$ it holds:
$$
f(a,b) = \left\{
           \begin{array}{ll}
              a, & if \; \; b= 1  \\
             b, & if \; \; a= 1  \\
             \end{array}
           \right.
$$
and otherwise,
$$ f(a,b) \geq   f(a,b-1) + f(a-1,b)+1.$$
\EL
\BPF
The proof is by induction on $a$ and $b$.
If $a = 1$ or $b = 1$ then the lemma follows directly from Claim~\ref{f_1}.
Otherwise, using the induction hypothesis, let $D'_{a,b-1}$ be a triangular submatrix of size  $f(a,b-1)$,
with row indices of size $a$ and column indices of size $b-1$,
and let $D''_{a-1,b}$ be a triangular submatrix of size $f(a-1,b)$,
with row indices of size $a-1$ and column indices of size $b$.

Assume that each row index of $D'_{a,b-1}$ is disjoint from all column indices of $D''_{a-1,b}$
(this is always possible for a large enough range of elements for the indices),
and let $x$ be a new element that does not appear in any of the row or column indices of $D'_{a,b-1}$ or of $D''_{a-1,b}$.
Add $x$ to each column index of $D'_{a,b-1}$ and to each row index of $D''_{a-1,b}$.
Therefore, the row and column indices of these two matrices are now subsets of size $a$ and $b$, respectively,
and each column index of $D'_{a,b-1}$ intersects all row indices of $D''_{a-1,b}$ (as they all contain $x$).

Now add to $D'_{a,b-1}$ one more row and column, as a last row and column,
defined by the row index  $\{x\} \cup S$, and the column index $\{x\} \cup T$,
where $S$ is a subset of size $a-1$ and $T$ is a subset of size $b-1$,
and $S$ and $T$ are disjoint from all row and column indices of $D'_{a,b-1}$.
Denote the resulting matrix  by $\tilde D'_{a,b-1}$.

Consider the following triangular matrix $D_{a,b}$ defined by all row and column indices of $\tilde D'_{a,b-1}$ and $D''_{a-1,b}$
(after adding $x$ and the additional row and column as described above):
$$
D_{a,b} = \left(
\begin{array}{cc}
 \tilde D'_{a,b-1} & O  \\
 J &D''_{a-1,b}
\end{array}
\right),
$$
where $J$ is the all-ones matrix and $O$ is the all-zeros matrix.
The size of $D_{a,b}$ is $f(a,b-1) + f(a-1,b)+1$,
and as stated, the row and column indices are subsets of size $a$ and $b$, respectively.
Hence, $f(a,b) \geq f(a,b-1) + f(a-1,b)+1$ as claimed.
\EPF\\

Before we solve this recursion, we need to recall the following definitions about recursion trees
that are useful to describe the expansion of a recursive formula.
A {\em rooted tree}, is a directed tree that has one node designated as the  root of the tree, and all edges are directed away from the root.
If $(u,v)$ is a directed edge in a directed tree, then $v$ is the child of $u$ in the tree. A {\em leaf} in the tree is a node with no edges coming out of it.
A node that is not a leaf is called an {\em internal} node of the tree.
A rooted tree is called a {\em full binary} tree if each node that is not a leaf  has exactly two children.
It is well known, and easy to prove by induction, that the number of leaves of a full binary tree is one more than the total number of internal nodes of the
tree.

Now, using these definitions and the recursion given in Lemma~\ref{const:triangular}, we can prove the following lower bound on $f(a,b)$.
From this bound follows immediately that $f(t,t) \geq {2t \choose t}-1$ as claimed.

\BL
\label{t_t}
For any $a,b \geq 1$ and a large enough $k$, $f(a,b) \ge  {a+b \choose a}-1$.
\EL
\BPF
If $b=1$ then by Claim~\ref{f_1} we have that $f(a,1) = a = {a+1 \choose a} -1$, and similarly if $a = 1$.
Therefore, assume that $a,b > 1$, and thus by Lemma~\ref{const:triangular}, $f(a,b) \geq f(a,b-1) + f(a-1,b)+1$.
The solution of this recursion is similar to the following recursion defined by Pascal's identity:
$${a + b \choose a} = {a + b-1 \choose a} + {a+ b-1 \choose a-1} = {a + b-1 \choose b-1} + {a+ b-1 \choose a-1}.$$
The only difference is the stopping conditions and the fact that the recursive formula for $f(a,b)$ has a plus one term.
Therefore, if we want to solve the recursion for $f(a,b)$,
we can expand instead the recursion for ${a+b  \choose a}$, and take into account the differences.

Let $T$ be a rooted binary labeled tree describing the expansion of the recursion $f(a,b) \geq f(a,b-1) + f(a-1,b)+1$,
where  $f(a,b)$ is the label of the root of the tree, and $f(a,b-1)$ and $f(a-1,b)$ are the labels of the two children of $f(a,b)$.
In general the children of a node labeled by $f(p,q)$ will be labeled by $f(p,q-1)$ and $f(p-1,q)$.
The labels of the leaves of the tree will be the stopping conditions of the recursion.

A similar tree $T'$, with the same structure as $T$, can be used to describe the expansion of ${a + b \choose a}$ using Pascal's identity,
where the labels are the binomial coefficients expanded by the recursion.
Since $T'$ describes the expansion of ${a + b \choose a}$, then the sum of the labels of its leaves is exactly ${a + b \choose a}$.

Now in order to solve the recursion for $f(a,b)$,
note that for each stopping term we loose $1$ compared to the expansion of Pascal's identity, since
$f(a,1) = f(1,a) = a$, whereas ${a+1 \choose 1} = {a+1 \choose a} = a+1$.
Thus, we have to subtract $1$ for each leaf of the tree $T$ from the sum ${a+b \choose a}$, for a total of $\ell$ ones,
where $\ell$ is the number of leaves of $T$.

However, the recursion for $f(a,b)$ has a plus $1$ term in each step of the recursion that is not a stopping condition,
whereas the recursion of ${a+b \choose a}$ does not have such a term.
Thus, we should sum these ones and add them to the total summed in the leaves.
The number of such ones that we should add is equal to the number of internal nodes of $T$, since each internal node corresponds
to a recursive step. But $T$ is a full binary tree, and thus, it has $\ell -1$ internal nodes.
Summarizing the above discussion, we get that:
$$
f(a,b) \geq {a+b \choose a} -\ell + \ell-1 = { a+b \choose a}-1.
$$
\EPF\\

\section{Constructions of large isolation sets  for $k \geq 2t$}
\label{sec:constructions}

In this section we prove Theorem~\ref{theorem:constructions},
and give  constructions of families of large isolation sets in $A_{k,t}$,
where for a large enough $k$, the constructions are the best possible, as we get an isolation set of size $k$.
The proof of the theorem contains several parts, according to the range of values of $k$ compared to $t$.
We first provide a basic construction of isolations sets of size $k-t+1$ for $k \geq 3t-2$,
and then use this construction to build large isolations sets
for  $2t \leq k \leq 3t-3$, for $3t-2 \leq k \leq 4t-3$, and finally for $k \geq 4t-3$.

\subsection{A construction of isolation sets of size $k-t+1$ for $k \geq 3t-2$}

We now prove that if $k \geq 3t-2$ then there exists an isolation set of size $k-t+1$ in $A_{k,t}$.
We first need to show that there exists an isolation matrix, not necessarily in $A_{k,t}$, of a certain structure,
such that each row and column of this matrix has the same number of ones.

\BCM
\label{claim:fool2t}
For any $q \geq p-1$, there exists an isolation matrix $F_{p,q}$ of size $(p + q)\times(p+ q)$,
such that there are $p$ ones and $q$ zeros in each column of $F_{p,q}$.
\ECM
\BPF
Take the circulant matrix $F_{p,q}$, whose first column is $(\overbrace{1,1,\cdots,1}^p,\overbrace{0,0,\cdots,0}^q)$.
It is not hard to verify that $F_{p,q}$ is an isolation matrix when $q \geq p-1$ (when $q = p-1$ the matrix is skew-symmetric).
Also, each column of $F_{p,q}$ is a cyclic permutation of the first column, and thus, each column contains $p$ ones and $q$ zeros.
See for example Figure~\ref{fig:fool2t}.
\EPF\\

\begin{figure}[htb!]
\captionsetup{width=0.8\textwidth}
	$$
	F_{5,4} =
\left(
\begin{array}{ccccccccc}
 {\bf 1} & 0 & 0 & 0 & 0  & 1 & 1 & 1 & 1   \\
 1 & {\bf 1} & 0 & 0  & 0 & 0 & 1 & 1& 1  \\
 1 & 1 & {\bf 1} & 0  & 0 & 0 & 0 & 1 & 1   \\
 1 & 1 & 1 & {\bf 1}  & 0 & 0 & 0 & 0 & 1   \\
 1 & 1 & 1 & 1 & {\bf 1}  & 0 & 0 & 0 & 0   \\
 0 & 1 & 1 & 1 & 1  & {\bf 1} & 0 & 0 & 0 \\
 0 & 0 & 1 & 1 & 1  & 1 & {\bf 1} & 0 & 0  \\
 0 & 0 & 0 & 1 & 1  & 1 & 1 & {\bf 1} & 0    \\
 0 & 0 & 0 & 0 & 1  & 1 & 1 & 1 & {\bf 1}   \\
\end{array}
\right)
	$$
\caption{An isolation matrix $F_{5,4}$ of size $(p + q)\times (p+q)  = 9\times 9$,
with $p = 5$ ones  and $q = 4$ zeros in each column.
The isolation set contains the ones in bold on the diagonal of $F_{5,4}$.}
\label{fig:fool2t}	
\end{figure}

\BL
\label{Lemma:isolation 3t-2}
If $k \geq 3t-2$ and $t\geq 2$, there exists an isolation set of size $k-t+1$ in $A_{k,t}$.
\EL
\BPF
Let $F_{p,q}$ be the isolation matrix described in Claim~\ref{claim:fool2t}, with $p = t$ and
$q = k-2t +1 \geq 3t-2-2t+1 = t-1 = p-1$. Let $I_{q+p}$ be the identity matrix of size $(q + p)\times (q + p)$,
$J_{q+p,p-1}$ the all-ones matrix of size $(q+p) \times (p-1)$, and $O_{p-1,q+p}$ the all-zeros matrix of size $ (p-1)\times (q +p)$.
Finally, let $X$ and $Y$ be the following matrices achieved by concatenating the above matrices as follows:
$$X = [I_{q+p} J_{q + p,p-1}], \; \; \; \; \; \;
Y = \left[
  \begin{array}{c}
    F_{p,q} \\
    0_{p-1,q+p} \\
  \end{array}
\right].
$$
Observe that $X Y = F_{p,q}$.
Furthermore, since each row of  $X$ and each column of $Y$ are vectors of length $q+2p-1 = k$ with exactly $p = t$ ones,
then we can view them as the characteristic vectors of subsets in ${[k] \choose t}$.
Thus, $XY = F_{p,q}$ is an isolation submatrix of $A_{k,t}$ of size $(q + p) \times (q + p) = (k-t+1)\times (k-t+1)$ as required.
\EPF\\

\subsection{A construction of large isolation sets  for $2t \leq k \leq 3t-3$ }

\BL
\label{isolation:small k}
Let $t\geq 2$ and $k = 2t+r$, where $0 \leq r \leq t-3$.
There exists an isolation matrix in $A_{k,t}$ of size $(2r+3)\times (2r+3)$.
\EL
\BPF
Let $t' = r+2$ and $k' = 2t'+r$. Thus, $k' = 2t'+r = 3t'-2$, and therefore, by
Lemma~\ref{Lemma:isolation 3t-2}, there exists an isolation matrix $F'$ of size $(k'-t'+1)\times (k'-t'+1)$
in $A_{k',t'}$, where the row and column indices of $F'$ are subsets of size $t'$ of $[k']$.

Since $k -k' = 2t+r - 2t'-r = 2 (t-r-2)$, there are still $2(t-r-2)$ elements from $[k]$ that were not used
to construct the row and column indices of $F'$.
Add to each row index of $F'$ half of these elements, and to each column index the other half.

Now the row and column indices are subsets of $[k]$ of size $t' + t-r-2 = t$,
and the resulting matrix is an isolation matrix of size $(2r+3)\times(2r+3)$ in $A_{k,t}$, as $2r+3 = k'-t'+1$.
\EPF

\subsection{A construction of large isolation sets  for $3t-2 \leq k \leq 4t-3$ }

\BL
\label{isolation:bigk}
Let $t\geq 2$ and $k = 2t+r$, where $t-2 \leq r \leq 2t-3$.
There exists an isolation matrix in $A_{k,t}$ of size $(2r+3)\times (2r+3)$.
\EL
\BPF
If $k = 3t -2$ then by Lemma~\ref{Lemma:isolation 3t-2}, there exists an isolation matrix of size $(2r+3)\times (2r+3)$ as required,
since $k-t+1 = 2t-1 = 2r+3$.
Otherwise,  $k > 3t-2$ and $r > t-2$ and define $k' = 3t-2$. Also let $O$ be the all-zero matrix, $J$  the all-one matrix,
and $F'$  the isolation submatrix of $A_{k',t}$, of size  $(k'-t+1)\times (k'-t+1) = (2t-1)\times (2t-1)$,
as promised by Lemma~\ref{Lemma:isolation 3t-2}.
Finally,  let $F''$ be another isolation matrix of size $(2r-2t+4)\times (2r-2t+4)$ that has the following structure:

$$
F'' =
 \left(
  \begin{array}{cccc|cccc}
    1 & 1 & \cdots & 1    & 1 & 0 & \cdots & 0 \\
    0 & 1 & \cdots & 1    & 1 & 1 &0 & 0 \\
    \vdots & 0 & \ddots & \vdots    & \vdots & \vdots & \ddots & 0 \\
    0 & 0 & 0 & 1    & 1& 1 & \cdots & 1 \\\hline
    0 & 0 &\cdots & 0  & 1 & 1 & \cdots & 1 \\
    1 & 0 & \cdots & 0  & 0 & 1 &\cdots & 1 \\
    \vdots & 1 & \ddots & \vdots  & \vdots & 0 & \ddots & \vdots \\
   1 & \cdots& 1 & 0  & 0 & \cdots & 0 & 1 \\
  \end{array}
\right)
$$
\\

We now show how to construct an isolation matrix $F$ of size $(2r+3)\times (2r+3)$ of the following structure
(the dimensions of the submatrices of $F$ are specified alongside the figure):

\begin{frame}
\footnotesize
\setlength{\arraycolsep}{3pt} 
$$
\begin{array}{rcl}
&
           \begin{array}{ccccccccccccccc}&&&&&&&&&&&&&&\end{array}
           \begin{array}{cc}
                      \overbrace{\begin{array}{cccccc}&&&&&\end{array}}^{r-t+2}
                      &
                      \overbrace{\begin{array}{cccccc}&&&&&\end{array}}^{r-t+2}
           \end{array}
&
\\
\begin{array}{r}
		\\
             \begin{array}{c}\\\\\end{array} \\
             \begin{array}{c}\\\\\\\end{array} \\
		r-t+2\left \{ \begin{array}{c}\\\\\\\\\end{array} \right.\\
             r-t+2\left \{ \begin{array}{c}\\\\\\\\\end{array} \right.
\end{array}
&
\left(
   \begin{array}{c|c}
               \begin{array}{ccccccccc} \\\\&&&&F' &&&&\\\\\end{array}
               &
               \begin{array}{c|c}
                           \begin{array}{ccccc}\\\\&&J &&\\\\\end{array}
                            &
              		  \begin{array}{c}
                          		     \begin{array}{ccccc}\\&&J &&\\\\\end{array}\\
				   		      \hline
                                            \begin{array}{ccccc}\\&&O &&\\\\\end{array}
                            \end{array}
                \end{array}
   \\
   \hline
   \begin{array}{c}
		 \begin{array}{ccccccccc} \\\\&&&&O &&&&\\\\\end{array}
		              \\ \hline
		              \begin{array}{c|c}
                                      \begin{array}{ccccc}\\&&O &&\\\\\\\end{array}&
				            \begin{array}{ccccc}\\&&J &&\\\\\\\end{array}
                           \end{array}
              \end{array}
   & F''
   \end{array}
\right)
&
\begin{array}{l}
             \left. \begin{array}{c}\\\\\\\end{array} \right\}t-1\\
             \left. \begin{array}{c}\\\\\\\end{array} \right\}t\\
		\begin{array}{c}\\\\\\\\\end{array}\\
             \begin{array}{c}\\\\\\\\\end{array}
\end{array}
\\
&
           \begin{array}{cc}
                      \underbrace{\begin{array}{cccccc}&&&&&\end{array}}_{t-1}
                      &
                      \underbrace{\begin{array}{cccccc}&&&&&\end{array}}_{t}
           \end{array}
           \begin{array}{ccccccccccccccc}&&&&&&&&&&&&&&\end{array}
&
\end{array}
$$
\end{frame}
\\

Since the sum of dimensions of $F'$ and $F''$ is $(2t-1) + (2r-2t+4) = 2r+3$, then $F$ is a matrix of size $(2r+3)\times (2r+3)$ as claimed.
In what follows we show that there is a way to assign row and column indices that are all subsets of ${[k] \choose t}$,
such that we get the above structure of $F',F''$ and $F$.
Then we can conclude that $F$ is an isolation submatrix of $A_{k,t}$,
since this structure of $F', F''$ and $F$, guaranties that any two ones on the diagonal of $F$ are not
in an all-ones submatrix of size $2\times 2$. \\

{\bf The row and column indices of $F'$:}
Denote the row and column indices of $F'$ by $R_1,...,R_{2t-1}$  and $C_1,...,C_{2t-1}$, respectively.
According to the construction described in Lemma~\ref{Lemma:isolation 3t-2}, both the row and column indices of $F'$ are
subsets of size $t$ of $[k']$  defined as follows:
\begin{itemize}
\item
For $i = 1,...,2t-1$: $R_i = \{i\} \cup S'$, where  $S' = \{2t,2t+1,...,3t-2\}$.
\item
For $i = 0,...,2t-2$,
$$C_{i+1} = \{i \bmod(2t-1)+1,\;(i+1) \bmod(2t-1)+1,\;...,(i+t-1) \bmod (2t-1)+1\}.$$
\end{itemize}
Note that the largest element in a column index of $F'$ is $2t-1$.
Furthermore, it appears in exactly the last $t$ column indices of $F'$.\\

{\bf The row and column indices of $F''$:}
Let $r' =  r-t+2$ and denote the row and column indices of $F''$ by $R_{2t},...,R_{2t+2r'-1}$ and  $C_{2t},...,C_{2t+2r'-1}$,
where:
\begin{itemize}
\item For $i = 0,...,2r'-1$: $C_{2t+i} = \{k-2r'+1 + i\} \cup S'',$ where $S'' = \{1,2,...,t-1\}$.
\item
For $i = 0,...,r'-1$:
$$R_{2t+i}  =  \{k-2r'+1+i,k-2r'+2+i,...,k-2r'+ r' +1 +i\}  \cup T,$$
where $T = \emptyset$ if  $r = 2t-3$, and otherwise, $T = \{2t,2t+1...,4t-r-4\} \subset S'$.
Note that the indices are well defined as $k-2r'+1 = 4t-r-3$, and the maximal element in $T$ is $4t-r-4$.
Furthermore, each index is a subset of size $r'+1 + |T| = r-t+3 + (2t-r-3) = t$ as required.
\item
$R_{2t+r'} = (R_{2t+r'-1} \setminus \{k-r'\})\cup \{2t-1\}.$
\item
For $i = 1,...,r'-1$: $ R_{2t+r'+i} =  (R_{2t+r'+i-1}\setminus\{k-r'+i\}) \cup \{k-2r'+i\}.$
\end{itemize}
It is not hard to verify that $F''$ has the structure described above, and that all row and column indices are
subsets of ${[k] \choose t}$. Therefore, $F''$ is an isolation submatrix of $A_{k,t}$ of size $(2r')\times (2r')$ as required.

Now if we consider the matrix defined by all the row and column indices  $R_1,...,R_{2t+2r'-1}$ and $C_1,...,C_{2t+2r'-1}$,
then we get the matrix $F$ as above.
To verify that $F$ has the structure claimed,
note that the first $r' = r-t+2$ row indices of $F''$, that is, $R_{2t},...,R_{2t+r'-1}$,  do not intersect with any of the column indices of $F'$,
since the largest element in a column index of $F'$ is $2t-1$, and the smallest element in these row indices is $x = \min \{2t, k-2r'+1\} = 2t$,
as $r \leq 2t-3$, and so $k-2r'+1 = k - 2(r-t+2) + 1 = 4t-r-3 \geq 2t$.

As to the row indices $R_{2t+r'},...,R_{2t+2r'-1}$, they intersect the last $t$ column indices of $F'$,
whereas, the column indices $C_{2t+r'},...,C_{2t+2r'-1}$ intersect with row indices $R_1,...,R_{t-1}$ of $F'$.
See also Figure~\ref{fig:foolt=4} for an example.
\EPF\\

\begin{figure}[htb!]
\captionsetup{width=0.8\textwidth}
\small{
$$
F =
\begin{array}{|c|c|c|c|c|c|c|c||c|c||c|c|}
\hline
         & 4 & 5 & 6 & 7 & 1 & 2 &  3 & 3 &3 &3 &3 \\
         & 3 & 4 & 5 & 6 & 7 & 1 &  2 & 2 &2 & 2& 2   \\
         & 2 & 3 & 4 & 5 & 6 & 7 &  1 & 1 &1 & 1&  1  \\
         & 1 & 2 & 3 & 4 & 5 & 6  & 7  & 9&10 &11 &12    \\\hline\hline
10, 9,8,1 & {\bf 1} &0 & 0 & 0 & 1 & 1 & 1 &1 &1 &1 &1  \\\hline
10, 9,8,2  & 1 & {\bf 1} & 0 & 0 & 0 & 1 & 1 &1 &1 &1 &1  \\\hline
10, 9,8,3 & 1 & 1 & {\bf 1} &0  & 0 & 0 & 1 & 1&1 &1 &1  \\\hline
10, 9,8,4 & 1 & 1 & 1 & {\bf 1} &0  & 0 & 0 & 1&1 &0 &0  \\\hline
10, 9,8,5 & 0 & 1 & 1 & 1 & {\bf 1}  & 0 & 0 & 1& 1& 0&0  \\\hline
10, 9,8,6  & 0 & 0 & 1 & 1 & 1 & {\bf 1} & 0 & 1&1 &0 &0  \\\hline
10, 9,8,7 & 0 & 0 & 0 & 1 & 1 & 1 & {\bf 1} & 1 &1  &0 &0  \\\hline\hline
 11, 10, 8, 9 & 0 & 0 & 0 & 0 & 0 & 0 & 0 &{\bf 1} & 1 & 1 & 0   \\\hline
 12,11, 8,10  & 0 & 0 & 0 & 0 & 0 & 0&  0&0 &{\bf 1} & 1&1   \\\hline\hline
 7,12, 8,11  & 0 & 0 & 0 & 1 & 1 & 1 & 1 & 0 &0 &{\bf 1}  &1    \\\hline
 9, 7, 8,12 & 0 & 0 & 0 & 1 & 1 & 1 & 1&1 &0  & 0  &{\bf 1}   \\\hline
\end{array}
$$
}
\caption{An isolation matrix  of size $(2r+ 3)\times (2r+3) = 11\times 11$ in $A_{k,t}$, where $k = 12$, $t=4$ and $r =k- 2t = 4$.}
\label{fig:foolt=4}	
\end{figure}

\subsection{A construction of maximal isolation sets  for $ k \geq 4t-3$ }
\BL
\label{isolation:maximal}
Let $t\geq 2$ and $k \geq 4t-3$. There exists an isolation matrix in $A_{k,t}$ of size $k\times k$.
\EL
\BPF
Let $k' = 4t-3$ and let $F$ be an isolation matrix of size $k'\times k'$,
with row and column indices that are subsets of size $t$ of $[k']$ as defined in the proof of Lemma~\ref{isolation:bigk}.
Now add $k-k'$ rows and $k-k'$ columns to $F$ with the following  indices:
\begin{itemize}
\item
For $i =1,...,k-k'$, add the row indices $\{k'+i,2t-1,2t,2t+1,...,3t-3\}.$
\item
For $i =1,...,k-k'$, add the column indices $\{k'+i,1,2,...,t-1\}$.
\end{itemize}
The resulting matrix is an isolation matrix of size $k \times k$. See Figure~\ref{fig:foolt=3} for an example.
\EPF\\

\begin{figure}[htb!]
\captionsetup{width=0.8\textwidth}
\small{
$$
F =
\begin{array}{|c||c|c|c|c|c||c|c||c|c||c|c|}
\hline
      & 3 & 4 & 5 & 1 & 2 & 2 & 2 & 2 & 2 & 2 &2\\
      & 2 & 3 & 4 & 5 & 1 & 1 & 1 & 1 &1 & 1 & 1\\
      & 1 & 2 & 3 & 4 & 5 & 6 & 7 & 8 &9 & 10 & 11\\\hline\hline
7,6,1  & {\bf 1} & 0 & 0 & 1 & 1& 1 & 1 & 1  &1 & 1&1 \\\hline
  7,6,2 & 1 &{\bf 1} & 0 & 0 & 1 & 1 & 1  & 1 &1 &1 &1 \\\hline
 7,6,3   & 1 & 1 & {\bf 1} & 0 & 0 & 1 & 1 & 0 &0 & 0 & 0 \\\hline
   7,6,4 & 0 & 1 & 1 & {\bf 1} & 0 & 1 & 1 & 0 &0 & 0 & 0 \\\hline
   7,6,5 & 0 & 0 & 1 & 1 & {\bf 1} & 1 & 1 & 0 &0 & 0 &0 \\\hline\hline
  8,7,6 & 0 & 0 & 0 & 0 & 0 & {\bf 1} & 1 &1 &0  & 0 & 0\\\hline
  9,8,7 & 0 & 0 & 0 & 0 & 0 & 0 & {\bf 1} &1 &1 & 0 &0  \\\hline\hline
   5,9,8 & 0 & 0 & 1 & 1 & 1 & 0 & 0 &  {\bf 1} &1 & 0 & 0\\\hline
   6,5,9 & 0 & 0 & 1 & 1 & 1 & 1 &0 &0 &{\bf 1} & 0& 0\\\hline\hline
    6,5,10 & 0 & 0 & 1 & 1 & 1 & 1 &0 &0 & 0&{\bf 1} & 0\\\hline
       6,5,11 & 0 & 0 & 1 & 1 & 1 & 1 &0 &0 & 0 & 0 &{\bf 1}\\\hline
\end{array}
$$
}
\caption{A maximal isolation matrix of size $k\times k =11 \times 11$ in $A_{k,t}$, where $k=11, t=3$.}
\label{fig:foolt=3}	
\end{figure}

\section{Bounds on the maximal size of isolation sets}

As we saw, the constructions given in Section~\ref{sec:constructions} are maximal for any $t\geq 2$ and $k \geq 4t-3$,
since we get an isolation set of size $k$ in this case. Our construction is also maximal for $k = 2t$.
In this section we prove Theorem~\ref{theorem:upperbounds}, and show that it is also maximal for $k=2t+1$.
We first need the following claims.

\BCM
\label{claim:nozeros1}
Let $k = 2t+1$ and let $F$ be an isolation matrix in $A_{k,t}$.
Then $F$ cannot contain a submatrix of size $2\times 2$ that is the all-zero matrix.
\ECM
\BPF
Assume, by contradiction, that $F$ has a submatrix of size $2\times 2$ that is the all-zero matrix,
and assume that this submatrix is defined by row indices $x,y$ and column indices $z,w$.
Assume, without loss of generality, that $x = \{1,2,...,t\}$.
Since $ x \cap z = x \cap w = \emptyset$ and $z \neq w$, we must have that  $z \cup w = \{t+1,...,2t+1\}$.
But we have also that $ y \cap z = y \cap w = \emptyset$, and therefore,
$y = \{1,2,...,t\}$. Thus, $y = x$ and this is a contradiction.
\EPF\\

\BCM
\label{claim:nozeros2}
Let $k = 2t+1$ and let $F$ be an isolation matrix in $A_{k,t}$. Then every row and column of $F$ has at most three zeros.
\ECM
\BPF
Assume, by contradiction, that $F$ has a row with four zeros,
and assume, without loss of generality, that it is the first row and that the zeros are in positions $2,3,4,5$ of this row.
Consider the following submatrix $W$ of $F$ defined by the first five rows of $F$ and columns $2,3,4,5$ of $F$:
$$
W =  \left(
  \begin{array}{c|ccccc}
          &    y_2 & y_3 & y_4 & y_5 \\\hline
    x_1   &   0 & 0 & 0 & 0 \\
    x_2   &  {\bf 1} &  &  &  \\
    x_3   &    & {\bf 1} &  &  \\
  x_4 &  &  & {\bf 1} &  \\
  x_5 &   &  &  & {\bf 1} \\
  \end{array}
\right)
$$
Let $W'$ be the submatrix containing the last four rows of the submatrix $W$.
Note that $W'$ is an isolation matrix of size $4 \times 4$.
If $W'$ contains two zeros in one of its rows, then with the zeros in the first row of $W$, we get that $F$ contains
a submatrix of size $2\times 2$ that is the all-zero matrix, and we get a contradiction by Claim~\ref{claim:nozeros1}.
Thus, $W'$ contains at most one zero in each one of its rows.
But since $W'$ is an isolation matrix of size $4\times 4$, it must contain at least ${4 \choose 2} = 6$ zeros, and again we get a contradiction.
\EPF\\

\BCM
\label{claim:regular}
Let $G = (V_1,V_2, E)$ be a $3$-regular bipartite graph, where $|V_1| = |V_2| = 6$. Then $G$ contains a $4$-cycle.
\ECM
\BPF
Let $V_1 = \{x_1,...,x_6\}$, $V_2 = \{y_1,...,y_6\}$, and assume, without loss of generality,  that $(x_1,y_1),(x_1,y_2),(x_1,y_3)\in E$
and  $(y_1,x_2),(y_1,x_3)\in E$.
If one of $y_2$ or $y_3$ is a neighbor of one of $x_2$ or $x_3$, then we are done since we have a $4$-cycle
(for example, if $(x_2,y_2)\in E$, then $(x_1,y_1,x_2,y_2,x_1)$ is a $4$-cycle).

Thus, consider now the case that $y_2$ and $y_3$ are not neighbors of  $x_2$ and $x_3$.
Therefore, each of $x_2$ and $x_3$ has two neighbors from $y_4,y_5,y_6$, and so they have a common neighbor, say $y_4$.
Since they are both also neighbors of $y_1$, we get a $4$-cycle $(x_2,y_1,x_3,y_4,x_2)$.
\EPF\\

\BL
Let $k = 2t+1$, $t\geq 3$, and let $F$ be an isolation matrix in $A_{k,t}$. Then the size of $F$ is at most $5\times 5$.
\EL
\BPF
Assume, by contradiction, that there is an isolation matrix $F$ of size $6 \times 6$,
and denote the rows of $F$ by $X_1,...,X_6$ and the columns of $F$ by $Y_1,...,Y_6$.
Since $F$ is an isolation matrix, then $X_i \circ Y_i = e_i$ for every $1 \leq i \leq 6$,
where $e_i$ is the $i^{\rm th}$ standard basis vector, and $\circ$ is the Hadamard (entry-wise) product.

First notice that $F$ cannot have a column $Y_i$ with  five ones, since then $X_i$ must have four zeros (as $X_i \circ Y_i = e_i$),
and this is impossible by Claim~\ref{claim:nozeros2}. Therefore, every column  of $F$ has at most four ones.
A similar argument holds for the rows of $F$.
Furthermore, if there exists a row/column with two ones then it has four zeros and again we get a contradiction.
Thus, every row and column of $F$ has at least three ones and at most four ones, and at least two zeros and at most three zeros.
We thus, have the following two cases:

{\bf Case 1:} Every row and column in $F$ has three ones.
Let $G$ be the bipartite $3$-regular  graph whose adjacency matrix is the complement of $F$
(that is,  each zero in $F$ is an edge of the graph).
Then by Claim~\ref{claim:regular}, the graph $G$ has a $4$-cycle.
Thus, $F$ has a submatrix of size $2 \times 2$ that is all zeros, and we get a contradiction by Claim~\ref{claim:nozeros1}.

{\bf Case 2 :} There exists at least one column in $F$ with four ones.
Assume, without loss of generality, that it is  $Y_1$ and that $Y_1 = (1,1,1,1,0,0)$.
But, $X_1 \circ Y_1 = e_1$ and  by Claim~\ref{claim:nozeros2} every row of $F$  contains at most three zeros.
Thus, $X_1 = (1,0,0,0,1,1)$.

Now consider the structure of the submatrix $W$ of $F$
defined by rows $X_2,X_3,X_4$ and columns $Y_2,Y_3,Y_4$ of $F$.
Notice that $W$ is an isolation matrix of size $3 \times 3$.
First we claim that there cannot be two zeros in any of the rows of $W$
(otherwise, we will get a submatrix of size $2\times 2$ of zeros with the zeros in $X_1$).
Hence, each row  of $W$ has at most one zero.
Also there cannot be two zeros in any of the columns of $W$,
since then we will get a $2\times 2$ all ones submatrix on the diagonal of $W$, in contradiction to $W$ being an isolation matrix.
Thus, each one of the rows and columns of $W$ must contain at most one zero. But since $W$ is an isolation matrix of size $3 \times 3$
it should have at least ${3 \choose 2 } = 3$ zeros, and so each one of the rows and columns of $W$ must contain exactly
one zero and two ones.
Therefore,  without loss of generality,   $F$ has the following structure:

$$
F =
\left(
  \begin{array}{cccccc}
  {\bf 1}   & 0   & 0   &  0     & 1   & 1       \\
  1   & {\bf 1}   & 0   &  1     &    &     \\
  1   & 1   & {\bf 1}  &  0     &    &    \\
  1    & 0   & 1   &  {\bf 1}     &    &     \\
   0    &    &   &       & {\bf 1}   &     \\
  0    &    &    &      &    & {\bf 1}    \\
  \end{array}
\right)
$$

Similar considerations as those above, show that there cannot be two zeros in positions $2,3,4$ of $X_5$ or of $X_6$
(otherwise, there will be a submatrix of size $2\times 2$ of zeros with the first row of $F$), and there cannot be three ones in positions $2,3,4$ of $X_5$
(otherwise, we get that $Y_5  = (1,0,0,0,1,1)$ and therefore $Y_6 = (1,1,1,1,0,1)$, or otherwise we get a submatrix of size $2\times 2$ of zeros.
 But then $Y_6$ contains five ones and again we get a contradiction). A similar argument holds for $X_6$. Thus, $X_5$ and $X_6$ each
 must contain one zero and two ones in positions $2,3,4$. Hence, without loss of generality, $F$ is of the following form
 (where columns $Y_5$ and $Y_6$ were determined according to $X_5,X_6$, so that $X_5 \circ Y_5 = e_5$, $X_6 \circ Y_6 = e_6$,
 and we do not get a submatrix  of size $2\times 2$ that is all-zeros):

$$
F =
\left(
  \begin{array}{cccccc}
 {\bf 1}   & 0   & 0   &  0     & 1   & 1       \\
  1   & {\bf 1}   & 0   &  1     & 0   & 0    \\
 1   & 1   & {\bf 1}  &  0     & 0   & 1   \\
 1    & 0   & 1   &  {\bf 1}     & 1  & 0    \\
 0    & 1   & 1   &  0     & {\bf 1}   &    \\
 0    & 1   & 0   &  1     &   & {\bf 1}    \\
  \end{array}
\right)
$$

Now denote  the row and column indices of $F$ by $x_1,...,x_6$ and $y_1,...,y_6$, respectively,
where each index is a subset of size $t$ of $[k] = [2t+1]$, and assume, without loss of generality, that $x_1 = \{1,...,t\}$.
From the structure of  $F$ we can deduce the following about its row and column indices:
\begin{itemize}
\item
Since $x_1 \cap y_2 = x_1 \cap y_3 = x_1 \cap y_4 = \emptyset$,
then $y_2,y_3,y_4 \subset \{t+1,...,2t+1\}$, $|y_i \cap y_j| = t-1$ for $2 \leq i \neq j \leq 4$,  and $|y_2 \cap y_3 \cap y_4| = t-2$.
Let $S = y_2 \cap y_3 \cap y_4$ and assume, without loss of generality, that $S = \{t+1,t+2,...,2t-2\}$,
and  $y_2 = S \cup \{2t-1,2t\}, y_3 = S \cup \{2t-1,2t+1\},y_4 = S \cup \{2t,2t+1\}$.
\item
Since $x_2 \cap y_3 = \emptyset$, $x_2 \cap y_2 \neq\emptyset$, $y_2 \cap y_3 = S \cup\{2t-1\}$ and $y_2 = S \cup \{2t-1,2t\}$ then  $2t \in x_2$.
In a similar way, $2t-1 \in x_3$, $2t+1 \in x_4$, $2t-1 \in x_5$ and  $2t \in x_6$.
\item
Furthermore, since $x_5 \cap y_1 = x_6 \cap y_1 = \emptyset$, then there exists a subset $T$ of size $|T| = t-1$ such that
$T \subseteq x_5 \cap x_6$. Since $x_5 \cap y_4 = x_6 \cap y_3 = \emptyset$ then $T \cap y_4 = T \cap y_3 = \emptyset$.
Thus, $T \subseteq \{1,2,...,t\}$.
\end{itemize}

Finally, since $F$ is an isolation matrix then either $x_5 \cap y_6 = \emptyset$ or $x_6 \cap y_5 = \emptyset$.
Assume first that $y_6 \cap x_5 = \emptyset$. Thus, from the above discussion and our last assumption,
$F$ has the following structure and row and column indices:
$$
F = \left(
  \begin{array}{c||c|c|c|c|c|c}
          & y_1 & y_2 = & y_3 = &  y_{4}=  & y_5 & y_6  \\
          &  &  S,2t-1,2t &  S,2t-1,2t+1 &  S,2t,2t+1 & &    \\\hline\hline
    x_1 =1,...,t  & {\bf 1}   & 0   & 0   &  0     & 1   & 1       \\\hline
     2t \in x_2  & 1   & {\bf 1}   & 0   &  1     & 0   & 0    \\\hline
    2t-1 \in x_3   & 1   & 1   & {\bf 1}  &  0     & 0   & 1   \\\hline
  2t+1 \in x_{4} & 1    & 0   & 1   &  {\bf 1}     & 1  & 0    \\\hline
   x_5 = T,2t-1 & 0    & 1   & 1   &  0     & {\bf 1}   &  0  \\\hline
   x_6 =T,2t & 0    & 1   & 0   &  1     &   & {\bf 1}    \\
  \end{array}
\right)
$$
But then $ Q \cap y_6 = \emptyset$, where $Q = \{2t,2t+1,2t-1\}\cup T$, and this is a contradiction, since  then $y_6 \subseteq [k]\setminus Q$
and $|[k]\setminus Q| = t-1$.

In a similar way, if $x_6 \cap y_5 = \emptyset$, and since also $x_3 \cap y_5 = \emptyset$,  then
$(\{2t-1,2t\}\cup T) \cap y_5 = \emptyset$.
On the other hand, $x_5 \cap y_5 \neq\emptyset$ and $x_5 = T \cup \{2t-1\}$ and again we get a contradiction.
\EPF\\

\subsection*{Acknowledgments}
This research did not receive any specific grant from funding agencies in the public, commercial, or not-for-profit sectors.

\bibliographystyle{plain}
\bibliography{biclique}

\end{document}